\documentclass[11pt,reqno]{amsart}

\usepackage[a4paper,margin=1.1in]{geometry}
\usepackage[T1]{fontenc}
\usepackage{lmodern}
\usepackage{amsmath,amssymb}
\usepackage{enumitem}
\usepackage[protrusion=true,expansion=false]{microtype}
\usepackage{xcolor}
\usepackage[colorlinks=true,linkcolor=blue,citecolor=blue,urlcolor=blue]{hyperref}

\theoremstyle{plain}
\newtheorem{theorem}{Theorem}            
\newtheorem{lemma}[theorem]{Lemma}
\newtheorem{proposition}[theorem]{Proposition}
\newtheorem{corollary}[theorem]{Corollary}

\theoremstyle{definition}
\newtheorem{definition}[theorem]{Definition}
\theoremstyle{remark}
\newtheorem{remark}[theorem]{Remark}
\theoremstyle{remark}
\newtheorem{example}[theorem]{Example}
\newcommand{\R}{\mathbb{R}}
\newcommand{\Null}{\mathcal N}
\newcommand{\ds}{\,ds}

\title[Geodesic connectedness with a causal Killing field]{Geodesic connectedness of globally hyperbolic spacetimes with a causal Killing vector field}

\author[R. Bartolo]{Rossella Bartolo}
\address{Dipartimento di Meccanica, Matematica e Management, Politecnico di Bari, Bari, Italy}
\email{rossella.bartolo@poliba.it}

\author[E. Caponio]{Erasmo Caponio}
\address{Dipartimento di Meccanica, Matematica e Management, Politecnico di Bari, Bari, Italy}
\email{erasmo.caponio@poliba.it}

\date{}
\dedicatory{To Miguel S\'anchez, on the occasion of his sixtieth birthday,
with great esteem and affection.}

\begin{document}

\begin{abstract}
We characterize the existence of a geodesic joining two prescribed points of a globally hyperbolic spacetime $(M,g)$ admitting a complete causal Killing vector 	field $K$ and a complete smooth spacelike Cauchy hypersurface $S$. We prove that 		two points $p,q\in M$ are joined by a geodesic if and only if there exists an 	$H^1$ curve $\gamma$ from $p$ to $q$ such that $g(\dot\gamma,K)$ is constant almost everywhere.  The statement is intrinsic and  allows $K$ to be timelike 	at some points and lightlike at others. A global splitting induced by the flow 	of $K$ is used in the proof, but no global bounds on the metric coefficients are imposed. The proof relies on stationary perturbations of $(M,g)$, whose relevant compactness properties follow from global hyperbolicity of $g$, and on uniform estimates based on an integral admissibility condition and a Glaeser-type inequality. Two  examples illustrate geodesic connectedness and the obstruction detected by the characterization.
\end{abstract}

\maketitle

\section{Introduction and statement of the result}

Let $(M,g)$ be a globally hyperbolic spacetime and let $K$ be a complete causal Killing vector field, namely $K$ is everywhere non-zero and $g(K,K)\le0$. 
By a theorem of Bernal and S\'anchez, every globally hyperbolic spacetime admits smooth, spacelike Cauchy hypersurfaces $S$, and in fact a smooth splitting; see \cite{BernalSanchez03}. 
 In this paper we impose the extra metric assumption that the chosen hypersurface $S$ is complete with respect to the induced Riemannian metric. We choose the time orientation so that $K$ is future-directed. Our main result shows that, under these geometric hypotheses,  the only obstruction to joining two prescribed points by a geodesic
 is the absence of an $H^1$ connecting curve with constant Killing charge. This characterization includes the mixed case in which the causal character of $K$ varies and, unlike the result obtained in \cite{BC}, does not require a global controllability assumption on the distribution associated with the kernel of the one-form induced by $K$ on $S$.

The two extreme cases in which $K$ is everywhere timelike or everywhere lightlike are now well understood, although for rather different reasons. If $K$ is everywhere timelike we can use $S$ and the flow of $K$ to obtain a global isometric splitting $S\times \R$. In this case  the spacetime is called {\em standard stationary} and the indefinite energy
functional associated with $g$ can be reduced, by using the natural Killing constraint, to a
functional defined on the  path space  of $S$.  Under suitable assumptions on the metric coefficients in the splitting, this  reduced functional is smooth, bounded from below and satisfies the Palais--Smale condition, hence the existence of critical points can be deduced. In  a standard stationary splitting this approach  was  introduced in  \cite{GiaMas91}, extending earlier work on the static case \cite{BeFoGi91}.
A later intrinsic reformulation was given in \cite{GiaPic99}: given $p,q\in M$, instead of fixing a global splitting, the energy functional is restricted to the constrained path space 
\begin{equation}\label{constK}
\begin{aligned}
\mathcal C_K(p,q)=\bigl\{\gamma\in H^1(I,M):\;&\gamma(0)=p,\ \gamma(1)=q,\\
&g(\dot\gamma,K)\equiv C_\gamma\ \text{a.e.\ for some }C_\gamma\in\R\bigr\},
\end{aligned}
\end{equation}
where $I=[0,1]$. Along a curve $\gamma$, the function $g(\dot\gamma,K)$ is the Noether charge associated with the symmetry $K$; we call it the \emph{Killing charge} of $\gamma$. Thus $\mathcal C_K(p,q)$ consists of the $H^1$ curves from $p$ to $q$ with constant Killing charge.

Under a pseudo-coercivity condition on this functional, the existence of geodesics connecting two fixed points can be obtained. However, as clarified in \cite{CaFlSa2008}, in the complete timelike Killing case this
intrinsic theory is essentially absorbed by the standard stationary one. Indeed, pseudo-coercivity\footnote{We mention that pseudo-coercivity has also been employed in a broader variational framework in \cite{CapCor23}, where the action functional of an autonomous indefinite Lagrangian admitting an infinitesimal symmetry with affine Noether charge is studied. A related variational framework for solutions with prescribed Lagrangian energy, connected with Fermat-type principles, has also been developed in \cite{CapCorGiaPic24}.} 
is closely tied to global hyperbolicity, and global hyperbolicity together with a complete timelike Killing field and
a complete spacelike Cauchy hypersurface yields a standard stationary
splitting  and the needed assumptions in the variational setup.

If, instead, $K$ is lightlike everywhere, the usual reduced energy functional
associated with a global isometric splitting induced by the flow of $K$ is no longer available. 
In particular, the Killing constraint no longer determines the time derivative of a curve satisfying it.
Nevertheless, in \cite[Theorem~1.2]{BaCaFl17} the existence of a connecting geodesic is still characterized by
the existence  of a $C^1$ curve $\gamma$
joining the prescribed endpoints and satisfying $g(\dot\gamma,K)\equiv C_\gamma$.
The proof in \cite{BaCaFl17} approximates the metric in a global splitting by stationary metrics and uses uniform estimates to obtain a limiting geodesic.
For related results on causality and geodesic connectedness in spacetimes admitting a causal Killing vector field, see \cite{AtGuHo20}.

Let us now recall why the variational reduction used in the standard stationary case
does not carry over directly to the causal case. 
Let $\Phi_t$ denote the flow of $K$. Its orbits are inextendible causal curves,
so each of them meets the spacelike Cauchy hypersurface $S$ exactly once.
Thus the map $\Phi:S\times\R\longrightarrow M$, $\Phi(x,t)=\Phi_t(x)$,
is a  diffeomorphism. Since the flow consists of isometries, each slice $\Phi(S\times\{t\})$ is a complete spacelike Cauchy hypersurface,
and, under this identification, the metric takes the form
\begin{equation*}
	g=g_0+\omega\otimes dt+dt\otimes\omega-\beta\,dt^2,
	\qquad \beta=-g(K,K)\ge0,
\end{equation*}
where $g_0$ is the Riemannian metric induced on $S$, with norm $|\cdot|_0$, and $\omega$ is a one-form on $S$. The Lorentzian character of $g$ is equivalent to
\begin{equation}\label{cono}
	\beta(x)+|\omega_x|_0^2>0,\qquad\hbox{for all }x\in S,
\end{equation}
see \cite[Proposition~3.3]{CaJaSa14}.
We shall call such a representation, obtained from the flow of the fixed complete Killing vector field $K$ and from a complete spacelike Cauchy hypersurface, a {\em complete Killing splitting}.
In this splitting $K=\partial_t$ and
\[
g(\nabla t,\nabla t)=-\frac{1}{\beta+|\omega|_0^2}<0.
\]
Consequently, $t$ is a temporal function, strictly increasing along every 
future-directed causal curve. The same observation applies to the stationary
perturbations introduced below. 
If $\beta>0$, the natural Killing
constraint in \eqref{constK} for a curve $\gamma(s)=(x(s), t(s))$,
\[
\omega(\dot x)-\beta(x)\dot t=C_\gamma,
\]
can be solved pointwise as
\[
\dot t=\frac{\omega(\dot x)-C_\gamma}{\beta(x)}.
\]
The endpoint condition on $t$ then determines $C_\gamma$, and the energy functional
restricted to the Killing-constrained path space can be written as a reduced
functional depending only on the spatial component $x$. In particular, the
reduced expression contains terms of the form
\[
\int_0^1 \frac{(\omega(\dot x)-C_\gamma)^2}{\beta(x)}\,ds
\]
and also the factor $\left(\int_0^1 \beta(x)^{-1}\,ds\right)^{-1}$.
Thus, the reduction is intrinsically connected to the strict timelike condition
$g(K,K)<0$, or equivalently $\beta>0$.

When $\beta$ is allowed to vanish, this procedure is no longer possible.
On the null locus
\begin{equation*}
\Null=\{x\in S:\beta(x)=0\},
\end{equation*}
the equation imposed by the Killing constraint ceases to determine $\dot t$,
becoming
\[
\omega(\dot x)=C_\gamma.
\]
Consequently, a spatial curve crossing $\Null$ or lying in  it can be lifted to a
Killing-constrained curve only if this condition holds on $\Null$.

The mixed causal case, where $\beta$ is non-negative and may vanish on a subset of $S$, was studied in \cite{BC}; see also the correction \cite{BC24}. There the proof was based on the approximation by standard stationary metrics as in \cite{BaCaFl17} and on affine control systems associated with the distribution $\mathcal D:=\ker\omega$. Together with global metric coefficient bounds, the result in \cite{BC} required $|\omega|_0$ to be bounded away from zero, $\mathcal D$ to be nowhere integrable, and $g_0(\nabla\beta,\omega^\sharp)=0$,
where $\omega^\sharp$ and $\nabla\beta$ are the vector fields $g_0$-dual to $\omega$ and $d\beta$, and $\nabla$ is the Levi-Civita connection of $g_0$.  
The role of the non-integrability assumption was to produce, for arbitrary spatial endpoints, suitable comparison curves. The orthogonality condition, on the other hand, was used in the estimates involving the derivatives $\dot t_\varepsilon$ of the time components of the approximating geodesics, where it removes the quadratic term (see \eqref{tddot} below).

The present paper keeps the same perturbative scheme, but changes the hypothesis which supplies the comparison curves. Indeed, we do not require a global controllability property of $\mathcal D$; instead, for the two prescribed endpoints we assume directly that the intrinsic path space $\mathcal C_K(p,q)$ is non-empty.
Moreover, the orthogonality assumption $g_0(\nabla\beta,\omega^\sharp)=0$ is removed thanks to a Glaeser-type estimate \cite{Glaeser63} for non-negative smooth functions.

As explained above, the usual stationary reduction does not apply directly when $\beta$ vanishes. Therefore we consider the metrics
\begin{equation}\label{geps}
 g_\varepsilon=g_0+\omega\otimes dt+dt\otimes\omega-\beta_\varepsilon dt^2,
 \qquad \varepsilon>0,
\end{equation}
with $\beta_\varepsilon:=\beta+\varepsilon$, 
for which $K=\partial_t$ is timelike and Killing, so that   $(S\times\R, g_\varepsilon)$ is a standard stationary spacetime for all $\varepsilon>0$. The goal is then to obtain, for sufficiently small $\varepsilon>0$, geodesics of $g_\varepsilon$ joining $p$ and $q$, together with estimates independent of $\varepsilon$, in order to pass to the limit as $\varepsilon\to 0$ in the geodesic equations of $g_\varepsilon$. The key compactness property comes from global hyperbolicity of the original metric $g$, rather than from global growth bounds on $\beta$ and $\omega$.
For a fixed point $q$ and a slice $S\times\{t_0\}$, this yields compactness of the intersection of the slice with the $g_\varepsilon$-causal past of $q$ for all sufficiently small $\varepsilon>0$; we do not assert global hyperbolicity of all these perturbed
spacetimes.

\smallskip
The following theorem is the main result of the paper.

\begin{theorem}\label{main}
Let $(M,g)$ be a globally hyperbolic spacetime endowed with a complete causal Killing vector field $K$ and a  complete (smooth, spacelike) Cauchy hypersurface $S$. Then, for any  $p,q\in M$, the following statements are equivalent:
\begin{enumerate}[label=\textup{(\roman*)},ref=\textup{(\roman*)},leftmargin=1.2cm]
\item\label{main-i} $p$ and $q$ are joined by a geodesic of $(M,g)$;
\item\label{main-ii} $\mathcal C_K(p,q)\ne\emptyset$.
\end{enumerate}
\end{theorem}

The implication \ref{main-i}$\Rightarrow$\ref{main-ii} follows at once from the  conservation law along a geodesic: indeed,
if $\gamma\colon I\to M$ is a geodesic of $g$ with $\gamma(0)=p$ and $\gamma(1)=q$, $g(\dot\gamma,K)$ is constant along $\gamma$ since $K$ is Killing. Hence $\gamma\in\mathcal C_K(p,q)$ and $\mathcal C_K(p,q)\ne\emptyset$. The converse is proved in Section~\ref{foce}. 

The paper is organized as follows. In Section~\ref{stationary} we obtain compactness of the intersections of causal past sets with a fixed slice for small stationary perturbations (Lemma~\ref{scudo}), introduce the reduced functional $J_\varepsilon$, and prove that it admits a minimizer for all sufficiently small $\varepsilon>0$. In Section~\ref{apriori}, writing $p=(x_0,t_0)$ and $q=(x_1,t_1)$ in a complete Killing splitting, we introduce admissible pairs $(\sigma,c)$, formed by an $H^1$ curve $\sigma$ in $S$ from $x_0$ to $x_1$ and a constant $c\in\R$ such that $\omega(\dot\sigma)=c$ almost everywhere on $\sigma^{-1}(\Null)$ and $(\omega(\dot\sigma)-c)^2/\beta(\sigma)$ is integrable on the complement. The spatial component of a curve in $\mathcal C_K(p,q)$, together with its Killing charge, is such a pair. Any admissible pair gives an upper bound, independent of $\varepsilon$, for the minima of $J_\varepsilon$, from which we derive uniform spatial $H^1$ bounds. We also prove a version of Glaeser's inequality on compact subsets of a Riemannian manifold. In Section~\ref{foce} this inequality gives a uniform bound on $\dot t_\varepsilon$; together with the conservation of $g_\varepsilon(\dot\gamma_\varepsilon,\dot\gamma_\varepsilon)$ along the approximating geodesics $\gamma_\varepsilon$, this yields uniform $C^1$ estimates for $\gamma_\varepsilon$ and Theorem~\ref{main} follows by passing to the limit.
Section~\ref{consequences} collects remarks and consequences of the characterization in Theorem~\ref{main}: its formulation in terms of the reduced functionals (Remark~\ref{variational-criterion}), the version of condition~\ref{main-ii} for $C^1$ curves whose Killing charge is nowhere zero or identically zero (Corollary~\ref{sign-criterion}), the invariance of geodesic connectability under the flow of $K$ (Corollary~\ref{orbit-invariance}), and simple sufficient criteria in a given splitting (Corollary~\ref{sufficient-criteria}). In particular, the non-integrability assumption of \cite{BC} yields, by the Chow--Rashevskii theorem, a curve tangent to $\mathcal D$ between any two points of $S$, as already observed in \cite[Section~1]{BC}; hence, under the hypotheses of Theorem~\ref{main}, this non-integrability assumption suffices for geodesic connectedness (Corollary~\ref{nonintegrable}). Example~\ref{positive-example} gives a geodesically connected spacetime which is not covered by Corollary~\ref{nonintegrable}, since the one-form induced by $K$ satisfies $\omega\wedge d\omega=0$ in every complete Killing splitting. In the splitting exhibited in the example, the orthogonality condition in \cite{BC} also fails. Since $K$ is complete, condition~\ref{main-ii} always holds when $K$ is everywhere timelike, by Corollary~\ref{sufficient-criteria}-(a); by Corollary~\ref{nonintegrable}, it can fail only if, in every complete Killing splitting, $\omega\wedge d\omega$ vanishes somewhere, as in the lightlike example of \cite[Section~7(c)]{BaCaFl17}. Example~\ref{negative-example} gives an analogous example in which $K$ is timelike on a non-empty open set, and shows that the relation of being joined by a geodesic is not transitive.

\begin{remark}\label{Deltat}
Once a complete Killing splitting has been fixed for the variational argument, we shall write
\[
\Delta t=t_1-t_0,
\]
and, after interchanging the two endpoints if necessary, we shall assume $\Delta t\ge0$. This is not restrictive because  if $\mathcal C_K(p,q)\ne\emptyset$,  reversing a curve produces an element of $\mathcal C_K(q,p)$ whose Killing charge is the opposite of the original one, and the property of being joined by a geodesic is symmetric.
\end{remark}

\section{Stationary perturbations and the reduced functionals}\label{stationary}
For $\varepsilon>0$, let us consider the metric $g_\varepsilon$ in \eqref{geps}. 
Notice that $g\prec g_\varepsilon$ and $g_\varepsilon\prec g_{\varepsilon'}$ if $0<\varepsilon<\varepsilon'$, where, for Lorentzian metrics $g'$ and $g''$, $g'\prec g''$ means that every non-zero $g'$-causal vector is $g''$-timelike; indeed, $g_{\varepsilon'}(v,v)=g_\varepsilon(v,v)-(\varepsilon'-\varepsilon)\,dt(v)^2$, and $dt(v)\neq0$ for every non-zero vector $v$ which is causal for $g$ or for some $g_\varepsilon$, since $t$ is a temporal function for all these metrics.

The following lemma shows that, for sufficiently small stationary perturbations, the intersection of the causal past of a point with a fixed slice is compact.
\begin{lemma}\label{scudo}
Assume that $(S\times\R,g)$ is globally hyperbolic, that each slice $S\times\{t\}$ is a Cauchy hypersurface for $g$, and that $(S,g_0)$ is complete. Let $t_0\in \R$ and $q=(x_1,t_1)$. Then there exists $\varepsilon_*>0$ such that, for every $\varepsilon\in(0,\varepsilon_*]$,
\begin{equation}\label{diamondcompact}
 J^-_\varepsilon(q)\cap (S\times\{t_0\})
\end{equation}
is compact in $S\times\{t_0\}$, where $J^-_\varepsilon$ denotes the causal past computed with respect to $g_\varepsilon$.
\end{lemma}
\begin{proof}
We adapt the limit-curve argument in \cite[Proposition~5.1]{BaCaFl17}. For $\varepsilon>0$ set
\[
A_\varepsilon=\{z\in S:(z,t_0)\in J^-_\varepsilon(q)\},
\]
so that the set in \eqref{diamondcompact} is $A_\varepsilon\times\{t_0\}$, and we have to prove that $A_\varepsilon$ is compact for all sufficiently small $\varepsilon>0$. Since $t$ is strictly decreasing along past-directed causal curves, $A_\varepsilon=\emptyset$ if $\Delta t<0$ and $A_\varepsilon=\{x_1\}$ if $\Delta t=0$. Assume henceforth that $\Delta t>0$.
	
We first prove that $A_\eta$ is bounded in $(S,g_0)$ for some $\eta>0$. Otherwise, there exist $\varepsilon_n\downarrow0$ and $z_n\in A_{\varepsilon_n}$ such that $d_0(z_n,x_1)\to+\infty$, where $d_0$ is the distance induced by $g_0$. Choose a past-directed $g_{\varepsilon_n}$-causal curve from $q$ to $y_n=(z_n,t_0)$, and extend it to the past along the integral curve of $-\partial_t$ through $y_n$. Parametrize the resulting past-inextendible curve $\gamma_n$ by arc length of the complete Riemannian metric $h=g_0+dt^2$, starting at $q$. If $r_n$ is the parameter at $y_n$, then
\[
r_n\ge d_0(x_1,z_n)\to +\infty,
\qquad t(\gamma_n(r))\ge t_0\quad\text{for }0\le r\le r_n.
\]
Since $\varepsilon_n\downarrow0$, all the curves $\gamma_n$ are $g_{\varepsilon_1}$-causal. Reversing the time orientation, the limit-curve theorem \cite[Theorem~3.1(1)]{Minguzzi08} (see also \cite[Proposition~3.31]{BeEhEa96}) gives, after passing to a subsequence, a past-inextendible $g_{\varepsilon_1}$-causal curve $\gamma\colon[0,+\infty)\to M$ issuing from $q$ such that $\gamma_n\to\gamma$ $h$-uniformly on compact subsets of $[0,+\infty)$. Since the metrics $g_{\varepsilon_n}$ have decreasing cones wider than those of $g$ and converge to $g$, it follows from \cite[Corollary~2.9]{Minguzzi08} that $\gamma$ is a past-directed $g$-causal curve. Moreover,  it is past-inextendible also for $g$, since a $g$-causal extension would be $g_{\varepsilon_1}$-causal. Local uniform convergence also gives $t(\gamma(r))\ge t_0$ for every $r\ge0$. On the other hand, extending $\gamma$ to the future along the integral curve of $\partial_t$ through $q$ we obtain an inextendible causal curve. Since  $S\times\{t_0\}$ is a Cauchy hypersurface, there is $r_*>0$ with $t(\gamma(r_*))=t_0$. As $t$ is strictly decreasing along the past-directed causal curve $\gamma$, we get $t(\gamma(r))<t_0$ for every $r>r_*$, a contradiction. Hence such an $\eta$ exists.
	
By completeness of $g_0$ and the Hopf--Rinow theorem, the set $\mathcal K=\overline{A_\eta}$ is compact. Set $\varepsilon_*=\eta$ and fix $\varepsilon\in(0,\varepsilon_*]$. Since the cones of $g_\varepsilon$ are contained in those of $g_\eta$, we have $A_\varepsilon\subset A_\eta\subset\mathcal K$, so it suffices to prove that $A_\varepsilon$ is closed. Let $z_m\in A_\varepsilon$ converge to $z\in S$, and let $\sigma_m$ be a future-directed $g_\varepsilon$-causal curve from $(z_m,t_0)$ to $q$. If $(w,t)$ lies on $\sigma_m$, then $t\ge t_0$, and following the integral curve of $\partial_t$ from $(w,t_0)$ to $(w,t)$ and then $\sigma_m$ up to $q$ gives a $g_\eta$-causal curve from $(w,t_0)$ to $q$. Hence, $w\in A_\eta\subset\mathcal K$. Thus every $\sigma_m$ lies in the compact set $\mathcal K\times[t_0,t_1]$. Parametrizing $\sigma_m$ by $t$, we write $\sigma_m(t)=(x_m(t),t)$, and the $g_\varepsilon$-causality of $\sigma_m$ gives
\[
\Bigl|\frac{dx_m}{dt}\Bigr|_0\le|\omega_{x_m}|_0+\sqrt{|\omega_{x_m}|_0^2+\beta_\varepsilon(x_m)}\le C,
\]
with $C$ independent of $m$, by compactness of $\mathcal K$. Hence the $h$-lengths of the curves $\sigma_m$ are uniformly bounded, and the limit curve theorem \cite[Theorem~3.1(2), case (i)]{Minguzzi08}, applied to $g_\varepsilon$, gives a $g_\varepsilon$-causal curve from $(z,t_0)$ to $q$. Therefore $z\in A_\varepsilon$, so $A_\varepsilon$ is closed, and hence compact.
\end{proof}
\begin{remark}
The previous proof is the only point where global hyperbolicity of the original metric $g$ is used to control the stationary perturbations. Without any global upper bound on $\beta$ or growth assumption on $\omega$, it yields compactness of the sets \eqref{diamondcompact} for sufficiently small $\varepsilon$.
\end{remark}
Next we deal with the reduced functionals: for each fixed $p,q\in M$ and $\varepsilon>0$, the restriction of the energy functional of $g_\varepsilon$ to the submanifold of curves $\gamma$ from $p$ to $q$ along which $g_\varepsilon(\dot\gamma,K)$ is constant is equivalent to a functional $J_\varepsilon$, the reduced functional, defined on the $H^1$ curves in $S$ joining the points $x_0$ and $x_1$, where $p=(x_0,t_0)$ and $q=(x_1,t_1)$.

Let
\[
\Omega_{x_0,x_1}(S)=\{x\in H^1(I,S): x(0)=x_0,\ x(1)=x_1\}
\]
and set
\begin{equation*}
 \mathcal A_\varepsilon(x)=\int_0^1\frac{\omega(\dot x)}{\beta_\varepsilon(x)}\ds,
 \qquad
 \mathcal B_\varepsilon(x)=\int_0^1\frac{1}{\beta_\varepsilon(x)}\ds,\quad x\in\Omega_{x_0,x_1}(S).
\end{equation*}
For $x\in\Omega_{x_0,x_1}(S)$, consider the lift $\gamma(s)=(x(s),t(s))$ from $p$ to $q$ along which $g_\varepsilon(\dot\gamma,K)$ is constant; its constant value $C_\varepsilon$ is determined by the endpoint condition on $t$. From
\begin{equation}\label{Ct}
	\dot t=\frac{\omega(\dot x)-C_\varepsilon}{\beta_\varepsilon(x)},
\end{equation}
we obtain
\begin{equation*}
 C_\varepsilon=\frac{\mathcal A_\varepsilon(x)-\Delta t}{\mathcal B_\varepsilon(x)}.
 \end{equation*}
 
The reduced functional is then
\begin{equation}\label{azione}
 J_\varepsilon(x)=\frac12\int_0^1|\dot x|_0^2\ds+\frac12\int_0^1\frac{\omega^2(\dot x)}{\beta_\varepsilon(x)}\ds
 -\frac{\bigl(\mathcal A_\varepsilon(x)-\Delta t\bigr)^2}{2\mathcal B_\varepsilon(x)} .
\end{equation}
Equivalently,
\begin{equation}\label{reazione}
 J_\varepsilon(x)=\frac12\int_0^1|\dot x|_0^2\ds+\frac12\int_0^1\beta_\varepsilon(x)\dot t^2\ds+C_\varepsilon\Delta t .
\end{equation}
\begin{remark}\label{Cepsilonmin}
A further representation makes the comparison estimates in Lemma~\ref{uniform} below immediate, namely 
\begin{align}
J_\varepsilon(x)
&=\frac12\int_0^1|\dot x|_0^2\ds+\min_{c\in\R}\mathcal H_\varepsilon(x,c),\nonumber \\
\mathcal H_\varepsilon(x,c)
&:=c\Delta t+\frac12\int_0^1
\frac{(\omega(\dot x)-c)^2}{\beta_\varepsilon(x)}\ds.\label{mathcalH}
\end{align}
Indeed,
\[
\partial_c\mathcal H_\varepsilon
=\Delta t-\mathcal A_\varepsilon(x)+c\mathcal B_\varepsilon(x),
\qquad
\partial_c^2\mathcal H_\varepsilon=\mathcal B_\varepsilon(x)>0.
\]
Thus the unique minimizer is $c=C_\varepsilon$, and substitution gives
\eqref{azione}.
\end{remark}
For all sufficiently small $\varepsilon>0$, the functional $J_\varepsilon$ attains its infimum, and each minimizer lifts to a geodesic of $g_\varepsilon$ joining $p$ and $q$, according to the following result.

\begin{proposition}\label{exist-min}
Under the assumptions of Theorem~\ref{main}, there exists $\varepsilon_*>0$ such that, for every $\varepsilon\in(0,\varepsilon_*]$, the functional $J_\varepsilon$ attains its infimum on $\Omega_{x_0,x_1}(S)$. If $x_\varepsilon$ is a minimizer and $t_\varepsilon$ is reconstructed by \eqref{Ct}, then $\gamma_\varepsilon=(x_\varepsilon,t_\varepsilon)$ is a geodesic of $(S\times\R,g_\varepsilon)$ joining $p$ and $q$.
\end{proposition}
\begin{proof}
By Lemma~\ref{scudo}, there exists $\varepsilon_*>0$ such that
$J^-_\varepsilon(q)\cap(S\times\{t_0\})$ is compact for every
$\varepsilon\in(0,\varepsilon_*]$. Thus the standard stationary metric
$g_\varepsilon$ satisfies the hypotheses of \cite[Theorem~4.2]{BaCaFl17}.
The proof of that theorem gives spatial $H^1$ bounds on each sublevel of
$J_\varepsilon$, first for $C^1$ curves and then for $H^1$ curves by density
and continuity of $J_\varepsilon$. As explained in \cite{BaCaFl17}, below Remark~5.2, the reduced functional is bounded from below,
satisfies the Palais--Smale condition and has complete sublevels, hence
attains its infimum. Its critical points lift, by \eqref{Ct}, to geodesics
with the prescribed endpoints.
\end{proof}

\section{A priori estimates}\label{apriori}
To pass to the limit in the family $(\gamma_\varepsilon)_{\varepsilon\in (0,\varepsilon_*]}$ obtained in Proposition~\ref{exist-min}, we first need a uniform comparison estimate. Fix a complete Killing splitting $M\simeq S\times\R$ and endpoints $p=(x_0,t_0)$, $q=(x_1,t_1)\in M$. We shall use the following notion of admissibility.
\begin{definition}
Let $\sigma\in H^1(I,S)$ and $c\in\R$. We define
\begin{equation}\label{Qbeta}
Q_\beta(\sigma,c)=\int_0^1\mathfrak{q}(\sigma(s),\dot\sigma(s),c)\ds,
\end{equation}
where, for $x\in S$ and $v\in T_xS$,
\[
\mathfrak{q}(x,v,c)=
 \begin{cases}
	 \dfrac{(\omega_x(v)-c)^2}{\beta(x)}, & \beta(x)>0,\\[0.9em]
	 0, & \beta(x)=0\ \text{and}\ \omega_x(v)=c,\\[0.2em]
	 +\infty, & \beta(x)=0\ \text{and}\ \omega_x(v)\ne c.
	 \end{cases}
\]
We shall say that the pair $(\sigma,c)$ is {\em admissible for $x_0,x_1$}, or that $\sigma$ is {\em admissible for the constant $c$}, if $\sigma\in\Omega_{x_0,x_1}(S)$ and $Q_\beta(\sigma,c)<+\infty$. When the endpoints are clear from the context, we simply say that $(\sigma,c)$ is admissible.
\end{definition}

\begin{remark}\label{misura}
For $\sigma\in H^1(I,S)$ and $c\in\R$, the integrand in \eqref{Qbeta} is measurable, so that $Q_\beta(\sigma,c)\in[0,+\infty]$ is well defined. Indeed, $\omega(\dot\sigma)$ is measurable, since $\dot\sigma$ is, and $\beta(\sigma)$ is continuous. Hence the integrand is measurable on the open set $\{s:\beta(\sigma(s))>0\}$, where it equals $(\omega(\dot\sigma)-c)^2/\beta(\sigma)$; on the closed set $\{s:\beta(\sigma(s))=0\}$ it equals $+\infty$ on the measurable subset where $\omega(\dot\sigma)\ne c$ and $0$ elsewhere. 

Moreover, admissibility is symmetric under the interchange of the endpoints, provided the constant $c$ is replaced by its opposite. Indeed, we note that $\mathfrak q(x,-v,-c)=\mathfrak q(x,v,c)$, thus we have $Q_\beta(\sigma(1-\cdot),-c)=Q_\beta(\sigma,c)$, so $(\sigma,c)$ is admissible for $x_0,x_1$ if and only if $(\sigma(1-\cdot),-c)$ is admissible for $x_1,x_0$. Hence,  assuming  $\Delta t\ge0$  after  interchanging $p$ and $q$ if necessary (recall Remark~\ref{Deltat}) is compatible with the existence of an admissible pair. We also note that  the interchanging of the endpoints  also leaves $c\Delta t$, and hence the quantity $U_{\sigma,c}$ of Lemma~\ref{uniform}, unchanged.
\end{remark}
\begin{remark}\label{admissible_from_Ck}
Let $\gamma\in\mathcal C_K(p,q)$ and write $\gamma=(\sigma,t)$ in a complete Killing splitting. If $c$ is the Killing charge of $\gamma$, then
\begin{equation}\label{splitting-constant}
 \omega(\dot\sigma)-\beta(\sigma)\dot t=c
 \qquad\text{a.e.\ on }I.
\end{equation}
Consequently $(\sigma,c)$ is admissible for $x_0,x_1$. Indeed, on $\{s: \beta(\sigma(s))>0\}$ one has
\[
 \frac{(\omega(\dot\sigma)-c)^2}{\beta(\sigma)}=\beta(\sigma)\dot t^{\,2},
\]
which is integrable because $\beta$ is bounded along the compact image of $\sigma$ and $\dot t\in L^2$. On the other hand, \eqref{splitting-constant} gives $\omega(\dot\sigma)=c$ a.e.\ on $\{s: \beta(\sigma(s))=0\}$.
\end{remark}

\begin{remark}\label{no-lift}
An admissible pair $(\sigma,c)$ need not admit an $H^1$ lift $\gamma:I\to M$ such that  $g(\dot\gamma, K)\equiv c$. Indeed, such a lift would require
\[
\frac{\omega(\dot\sigma)-c}{\beta(\sigma)}\in
L^2\bigl(\{s:\beta(\sigma(s))>0\}\bigr),
\]
whereas admissibility only yields square integrability of
$(\omega(\dot\sigma)-c)/\sqrt{\beta(\sigma)}$ on this set.
For instance, take $S=\R$, $g_0=dx^2$,
$\omega=dx$, $\beta(x)=x^2$, and
\[
\sigma(s)=e^s-1,\qquad c=1,\qquad 0\le s\le1.
\]
Then $\omega(\dot\sigma)-c=\sigma$ and $Q_\beta(\sigma,1)=1$, but any
such lift would satisfy
\[
\dot t(s)=\frac{1}{e^s-1},\qquad s>0,
\]
which is not in $L^2\big([0,1]\big)$. The equivalence established in
Corollary~\ref{admissibility-equivalence} below concerns only the existence of an admissible
pair and the non-emptiness of $\mathcal C_K(p,q)$.
\end{remark}

\begin{lemma}\label{uniform}
Let $\sigma\in H^1(I,S)$ be admissible for $c\in\R$, and let $\varepsilon_*>0$ be as in Proposition~\ref{exist-min}.
Then, for every $\varepsilon\in (0,\varepsilon_*]$,
\begin{equation}\label{Usigma}
 J_\varepsilon(\sigma)\le U_{\sigma,c}:=\frac12\int_0^1|\dot\sigma|_0^2\ds+c\Delta t+\frac12Q_\beta(\sigma,c).
\end{equation}
\end{lemma}

\begin{proof}
From Remark~\ref{Cepsilonmin} we have
\[
J_\varepsilon(\sigma)\le
\frac12\int_0^1|\dot\sigma|_0^2\ds+c\Delta t
+\frac12\int_0^1\frac{(\omega(\dot\sigma)-c)^2}{\beta_\varepsilon(\sigma)}\ds.
\]
The last integral is at most $Q_\beta(\sigma,c)$ since on $\{s:\beta(\sigma(s))>0\}$, $\beta_\varepsilon\ge\beta$ holds, while on $\{s:\beta(\sigma(s))=0\}$
Remark~\ref{misura} gives $\omega(\dot\sigma)-c=0$ almost everywhere.
This proves \eqref{Usigma}.
\end{proof}

When $\omega(\dot\sigma)\equiv c$, we have $Q_\beta(\sigma,c)=0$ and \eqref{Usigma} reduces to $J_\varepsilon(\sigma)\le\frac12\int_0^1|\dot\sigma|_0^2\ds+c\Delta t$. This is the analogue of the estimate obtained in \cite{BC} by means of trajectories of a family of  affine control systems depending on $\varepsilon$ (see \cite[Proposition~11 and Remark~12]{BC}).

The quantity $U_{\sigma,c}$ has a simple meaning when $(\sigma,c)$ comes from a curve $\gamma=(\sigma,t)\in\mathcal C_K(p,q)$ with Killing charge $c$.  In this case,  from Remark~\ref{admissible_from_Ck} we have $Q_\beta(\sigma,c)=\int_0^1\beta(\sigma)\dot t^{\,2}\ds$, and \eqref{splitting-constant} gives $U_{\sigma,c}=\frac12\int_0^1g(\dot\gamma,\dot\gamma)\ds$, the $g$-energy of $\gamma$. In general, $\min_cU_{\sigma,c}$ is the limit of $J_\varepsilon(\sigma)$ as $\varepsilon\to0^+$; see Remark~\ref{variational-criterion}.

We have the following   monotonicity  property for the family of reduced functionals.
\begin{lemma}\label{monotone}
For every fixed $x\in\Omega_{x_0,x_1}(S)$, the map $\varepsilon\mapsto J_\varepsilon(x)$ is non-increasing. 
\end{lemma}
\begin{proof}
If $0<\varepsilon\le\eta$, then $\beta_\eta\ge\beta_\varepsilon$ and, recalling the definition of $\mathcal H_\varepsilon$ in \eqref{mathcalH}, we have
$\mathcal H_\eta(x,c)\le\mathcal H_\varepsilon(x,c)$ for every $c\in\R$; thus Remark~\ref{Cepsilonmin} gives $J_\eta(x)\le J_\varepsilon(x)$.
\end{proof}
Now we are ready to prove some  uniform bounds involving the family of minimizers of the reduced functionals in Proposition~\ref{exist-min}.
\begin{lemma}\label{uniform2}
Under the  assumptions of Theorem~\ref{main}, suppose that there
exists an admissible pair $(\sigma,c)$ for $x_0,x_1$ in the chosen complete
Killing splitting.  For every
$\varepsilon\in(0,\varepsilon_*]$, let
$\gamma_\varepsilon=(x_\varepsilon,t_\varepsilon)$ be as in
Proposition~\ref{exist-min}, with $x_\varepsilon$ a minimizer of
$J_\varepsilon$. Then the family $(x_\varepsilon)_{\varepsilon\in (0,\varepsilon_*]}$ is bounded in $H^1(I,S)$ and its images are contained in a compact subset of $S$. Moreover, the families
\begin{equation*}
\left( \int_0^1\beta_\varepsilon(x_\varepsilon)\dot t_\varepsilon^2\ds\right)_{\varepsilon\in (0,\varepsilon_*]}\quad\text{and } \quad \left(C_\varepsilon\right)_{\varepsilon\in (0,\varepsilon_*]}
\end{equation*}
are bounded, where
\begin{equation}\label{Ceps}
C_\varepsilon=\omega(\dot x_\varepsilon)-\beta_\varepsilon(x_\varepsilon)\dot t_\varepsilon.
\end{equation}
\end{lemma}

\begin{proof}
Fix the admissible pair $(\sigma,c)$ and let $U_{\sigma,c}$ be as in
Lemma~\ref{uniform}. Set $\eta=\varepsilon_*$.
As recalled in the proof of Proposition~\ref{exist-min}, each sublevel of the fixed functional
$J_\eta$ has bounded spatial $H^1$ energy. If $0<\varepsilon\le\eta$, then Lemma~\ref{monotone} and the minimality of $x_\varepsilon$ give
\begin{equation}\label{fixedsublevel}
 J_\eta(x_\varepsilon)\le J_\varepsilon(x_\varepsilon)
 \le J_\varepsilon(\sigma)
 \le U_{\sigma,c}.
\end{equation}
Thus the entire family $(x_\varepsilon)$, with $0<\varepsilon\le\eta$, is contained in one sublevel of $J_\eta$ and has bounded $H^1$ energy, uniformly with respect to $\varepsilon$. By the Cauchy--Schwarz inequality, the $g_0$-lengths of these curves are uniformly bounded; since all of them start at $x_0$, completeness of $g_0$ and the Hopf--Rinow theorem give a compact
set $\mathcal K\subset S$ containing all the images.
Now set
\[
 Y_\varepsilon=\int_0^1|\dot x_\varepsilon|_0^2\ds,
 \qquad
 Z_\varepsilon^2=\int_0^1\beta_\varepsilon(x_\varepsilon)\dot t_\varepsilon^2\ds .
\]
We already know that $(Y_\varepsilon)_{\varepsilon\in (0,\eta]}$ is bounded. Since $\mathcal K$ is compact, there are constants $b_{\mathcal K},w_{\mathcal K}\ge0$ such that
\[
 0\le\beta_\varepsilon(x)\le b_{\mathcal K},
 \qquad |\omega_x|_0\le w_{\mathcal K},
 \qquad \text{for all }x\in \mathcal K,\ 0<\varepsilon\le\eta .
\]
Integrating the conservation law \eqref{Ceps} and using the Cauchy--Schwarz inequality yields
\begin{equation}\label{localCbound}
 |C_\varepsilon|
 \le w_{\mathcal K}\int_0^1|\dot x_\varepsilon|_0\ds
 +b_{\mathcal K}^{1/2} Z_\varepsilon
 \le A_{\mathcal K}+b_{\mathcal K}^{1/2}Z_\varepsilon,
\end{equation}
for a suitable $A_{\mathcal K}$.
On the other hand, from \eqref{reazione} and \eqref{fixedsublevel},
\begin{equation}\label{localenergy}
 \frac12Y_\varepsilon+\frac12Z_\varepsilon^2+C_\varepsilon\Delta t
 \le U_{\sigma,c} .
\end{equation}
As $\Delta t\ge0$, \eqref{localCbound} and \eqref{localenergy} imply
\[
 \frac12Z_\varepsilon^2
 \le U_{\sigma,c}+\Delta t\,|C_\varepsilon|
 \le U_{\sigma,c}+\Delta t\,A_{\mathcal K}+\Delta t\,b_{\mathcal K}^{1/2}Z_\varepsilon .
\]
This quadratic inequality gives a bound for $(Z_\varepsilon)_{\varepsilon\in (0,\varepsilon_*]}$. The bound for $(C_\varepsilon)_{\varepsilon\in (0,\varepsilon_*]}$ follows from \eqref{localCbound}.
\end{proof}

The next lemma is a version, on compact subsets of a Riemannian manifold, of the elementary estimate usually referred to as Glaeser's inequality. It is one of the basic ingredients of Glaeser's theorem on square roots of non-negative differentiable functions \cite{Glaeser63}; Dieudonn\'e \cite{Dieudonne70} gave a shorter proof of that theorem, based precisely on the pointwise estimate $|Df|^2\le Cf$ near the zero set. 
We include a short proof in the form needed below. It uses only
non-negativity and local bounds for the first two covariant derivatives of
$\beta$.
The estimate will be used in Section~\ref{foce} on the compact set $\mathcal K$ obtained in Lemma~\ref{uniform2}. 
\begin{lemma}\label{radice}
Let $\mathcal K\subset S$ be compact and let $\beta\in C^\infty(S)$ satisfy $\beta\ge0$. Then there exists $a_{\mathcal K}\ge0$ such that
\begin{equation}\label{chiave}
 |\nabla\beta(x)|_0^2\le a_{\mathcal K}\beta(x),\qquad\text{for all }x\in \mathcal K.
\end{equation}
\end{lemma}

\begin{proof}
Choose a relatively compact open neighborhood $V$ of $\mathcal K$ and
$r_0>0$ such that every unit-speed geodesic of the metric $g_0$  starting in $\mathcal K$
is defined and remains in $V$ for $0\le r\le r_0$. Choose $m_0>0$ so that
\begin{equation}\label{grabetabound}
m_0\ge\sup_V|\nabla^2\beta|_0,
\qquad m_0r_0\ge\sup_{\mathcal K}|\nabla\beta|_0,
\end{equation}
where the Hessian norm is the operator norm. Fix $x\in\mathcal K$ and set
$\kappa=|\nabla\beta(x)|_0$. If $\kappa=0$, \eqref{chiave} is immediate. Otherwise,
let $\psi$ be the unit-speed geodesic with
$\psi(0)=x$ and $\dot\psi(0)=-\nabla\beta(x)/\kappa$.
By \eqref{grabetabound}, $\kappa/m_0\le r_0$. Taylor's formula and the Hessian bound give
\[
0\le\beta(\psi(r))\le\beta(x)-\kappa r+\frac{m_0}{2}r^2,
\qquad 0\le r\le \kappa/m_0.
\]
Taking $r=\kappa/m_0$ yields
$|\nabla\beta(x)|_0^2\le2m_0\beta(x)$ on all of $\mathcal K$.
\end{proof}

\section{Uniform \texorpdfstring{$C^1$}{C1} bounds and proof of Theorem~\ref{main}}\label{foce}
Throughout this section we assume that $\mathcal C_K(p,q)\ne\emptyset$. By Remark~\ref{admissible_from_Ck}, an admissible pair exists, so Lemma~\ref{uniform2} applies. The estimates below use the intrinsic  hypothesis  through the existence of this pair. Let us specify the notation used in this section. The symbol $\omega^\sharp$ denotes the $g_0$-dual of $\omega$, and $W$ the $g_0$-skew-symmetric endomorphism associated with $d\omega$, namely
$d\omega(v,w)=g_0\big(v,W(w)\big)$.
Constants denoted by $C$ may change from line to line, but are always independent of $\varepsilon$.

For all $\varepsilon\in (0,\varepsilon_*]$ let $\gamma_\varepsilon=(x_\varepsilon,t_\varepsilon)$ be as in Proposition~\ref{exist-min}, with $x_\varepsilon$ a minimizer of $J_\varepsilon$.
The geodesic equations for the approximating metrics $g_\varepsilon$ are
\begin{equation*}
\begin{cases}
 \nabla_{\dot x_\varepsilon}\dot x_\varepsilon-
 \dot t_\varepsilon W(\dot x_\varepsilon)+\omega^\sharp(x_\varepsilon)\ddot t_\varepsilon+
 \dfrac12\dot t_\varepsilon^2\nabla\beta(x_\varepsilon)=0,\\[0.8em]
 \omega(\dot x_\varepsilon)-\beta_\varepsilon(x_\varepsilon)\dot t_\varepsilon=C_\varepsilon .
\end{cases}
\end{equation*}
Taking the scalar product of the first equation with $\omega^\sharp$ and
differentiating the conservation law gives
\begin{align*}
\bigl(\beta_\varepsilon(x_\varepsilon)
      +|\omega^\sharp(x_\varepsilon)|_0^2\bigr)\ddot t_\varepsilon
=&g_0(\dot x_\varepsilon,\nabla_{\dot x_\varepsilon}\omega^\sharp)
+\dot t_\varepsilon g_0(W(\dot x_\varepsilon),\omega^\sharp)\\
&\qquad-\dot t_\varepsilon g_0(\nabla\beta,\dot x_\varepsilon)-\frac12\dot t_\varepsilon^2g_0(\nabla\beta,\omega^\sharp),
\end{align*}
where the tensor fields on the right-hand side are evaluated along
$x_\varepsilon$. Equivalently,
\begin{equation}\label{tddot}
 \ddot t_\varepsilon=\mathsf A_\varepsilon+\mathsf B_\varepsilon\dot t_\varepsilon+\mathsf Q_\varepsilon\dot t_\varepsilon^2,
\end{equation}
where, along $x_\varepsilon$,
\begin{align*}
 \mathsf A_\varepsilon
 &=\frac{g_0(\dot x_\varepsilon,\nabla_{\dot x_\varepsilon}\omega^\sharp)}
 {\beta_\varepsilon(x_\varepsilon)+|\omega^\sharp(x_\varepsilon)|_0^2},\\
 \mathsf B_\varepsilon&=\frac{g_0(W(\dot x_\varepsilon),\omega^\sharp)-g_0(\nabla\beta,\dot x_\varepsilon)}{\beta_\varepsilon(x_\varepsilon)+|\omega^\sharp(x_\varepsilon)|_0^2},\\
 \mathsf Q_\varepsilon&=-\frac{g_0(\nabla\beta,\omega^\sharp)}
 {2\bigl(\beta_\varepsilon(x_\varepsilon)+|\omega^\sharp(x_\varepsilon)|_0^2\bigr)}.
\end{align*}
By \eqref{cono}, the continuous function $\beta+|\omega|_0^2$ has a positive minimum $\mu_{\mathcal K}$ on the compact set $\mathcal K$ of Lemma~\ref{uniform2}. Hence, 
\[
 \beta_\varepsilon(x_\varepsilon)+|\omega_{x_\varepsilon}|_0^2\ge \mu_{\mathcal K}
\]
for all $\varepsilon\in (0,\varepsilon_*]$. By the smoothness of the coefficients on $\mathcal K$, and since $|g_0(\nabla\beta,\omega^\sharp)|\le|\omega|_0\sqrt{a_{\mathcal K}\,\beta}$ on $\mathcal K$ by Lemma~\ref{radice} and the Cauchy--Schwarz inequality, there exists $C\ge 0$ such that
\begin{equation}\label{triade}
 |\mathsf A_\varepsilon|\le C|\dot x_\varepsilon|_0^2,
 \qquad
 |\mathsf B_\varepsilon|\le C|\dot x_\varepsilon|_0,
 \qquad
 |\mathsf Q_\varepsilon|\le C\sqrt{\beta(x_\varepsilon)} .
\end{equation}
In the next two lemmas we obtain, uniformly for $0<\varepsilon\le\varepsilon_*$,
bounds on $\|\dot t_\varepsilon\|_{L^\infty}$ and $\|\dot x_\varepsilon\|_{L^\infty}$.
\begin{lemma}\label{boundot}
The family $(\dot t_\varepsilon)_{\varepsilon\in (0,\varepsilon_*]}$ is bounded in $L^\infty(I)$.
\end{lemma}

\begin{proof}
First, from \eqref{triade} and Lemma~\ref{uniform2},
\begin{equation}\label{AB}
	\|\mathsf A_\varepsilon\|_{L^1}\le C
	\quad\text{ and }\quad
	\|\mathsf B_\varepsilon\|_{L^1}\le C .
\end{equation}
Setting $y_\varepsilon=\dot t_\varepsilon$, again by \eqref{triade} and Lemma~\ref{uniform2} we have
\begin{equation}\label{Qy}
\int_0^1|\mathsf Q_\varepsilon y_\varepsilon|\ds
\le C\int_0^1\sqrt{\beta(x_\varepsilon)}|y_\varepsilon|\ds
\le C\left(\int_0^1\beta(x_\varepsilon)y_\varepsilon^2\ds\right)^{1/2}\le C.
\end{equation}
Thus $\mathsf R_\varepsilon=\mathsf B_\varepsilon+\mathsf Q_\varepsilon y_\varepsilon$ is bounded in $L^1(I)$, and by \eqref{tddot} the function $y_\varepsilon$ solves the linear equation
\[
\dot y_\varepsilon=\mathsf A_\varepsilon+\mathsf R_\varepsilon y_\varepsilon .
\]
Since $\int_0^1 y_\varepsilon\ds=\Delta t$, there exists $s_\varepsilon\in I$ such that $y_\varepsilon(s_\varepsilon)=\Delta t$, and the variation-of-constants formula gives, for every $s\in I$,
\[
y_\varepsilon(s)=\Delta t\,\exp\Bigl(\int_{s_\varepsilon}^s\mathsf R_\varepsilon\Bigr)
+\int_{s_\varepsilon}^s\exp\Bigl(\int_r^s\mathsf R_\varepsilon\Bigr)\mathsf A_\varepsilon(r)\,dr .
\]
Hence $\|y_\varepsilon\|_{L^\infty}\le\exp\bigl(\|\mathsf R_\varepsilon\|_{L^1}\bigr)\bigl(\Delta t+\|\mathsf A_\varepsilon\|_{L^1}\bigr)$, which is bounded uniformly in $\varepsilon$ by \eqref{AB} and \eqref{Qy}.
\end{proof}

\begin{lemma}\label{boundotics}
The family $(\dot x_\varepsilon)_{\varepsilon\in (0,\varepsilon_*]}$ is bounded in $L^\infty(I)$.
\end{lemma}

\begin{proof}
By \eqref{reazione} and \eqref{Ceps}, $J_\varepsilon(x_\varepsilon)=\frac12\int_0^1g_\varepsilon(\dot\gamma_\varepsilon,\dot\gamma_\varepsilon)\ds$. Since $\gamma_\varepsilon$ is a geodesic of $g_\varepsilon$ defined on $I$, the function $g_\varepsilon(\dot\gamma_\varepsilon,\dot\gamma_\varepsilon)$ is constant and equal to $2J_\varepsilon(x_\varepsilon)$, which is at most $2U_{\sigma,c}$ by \eqref{fixedsublevel}. Writing $\omega(\dot x_\varepsilon)=C_\varepsilon+\beta_\varepsilon(x_\varepsilon)\dot t_\varepsilon$ by \eqref{Ceps}, we obtain, on $I$,
\[
|\dot x_\varepsilon|_0^2=2J_\varepsilon(x_\varepsilon)-2C_\varepsilon\dot t_\varepsilon-\beta_\varepsilon(x_\varepsilon)\dot t_\varepsilon^{\,2}\le2U_{\sigma,c}+2|C_\varepsilon|\,|\dot t_\varepsilon|,
\]
and the right-hand side is bounded uniformly in $\varepsilon$ by Lemmas~\ref{uniform2} and~\ref{boundot}.
\end{proof}

\begin{proof}[Proof of \ref{main-ii}$\Rightarrow$\ref{main-i} in Theorem~\ref{main}]
For $\varepsilon\in (0,\varepsilon_*]$, let $x_\varepsilon$ be a minimizer of $J_\varepsilon$ given by Proposition~\ref{exist-min}, and let $t_\varepsilon$ be reconstructed by \eqref{Ct}. Then $\gamma_\varepsilon=(x_\varepsilon,t_\varepsilon)$
is a geodesic of $(S\times\R,g_\varepsilon)$ joining $p$ and $q$. By Lemma~\ref{uniform2}, the curves $x_\varepsilon$ lie in the compact set $\mathcal K$; since $t_\varepsilon(0)=t_0$, Lemma~\ref{boundot} shows that the curves $\gamma_\varepsilon$ lie in a fixed compact set $\mathcal K\times[t_0-C,t_0+C]$. Moreover, $|\dot\gamma_\varepsilon|_h^2=|\dot x_\varepsilon|_0^2+\dot t_\varepsilon^{\,2}$, where $h=g_0+dt^2$ is the Riemannian metric introduced in the proof of Lemma~\ref{scudo}, so by Lemmas~\ref{boundot} and~\ref{boundotics} the curves $\gamma_\varepsilon$ are uniformly $h$-Lipschitz. Choose any sequence $\varepsilon_n\downarrow0$. By the Ascoli--Arzel\`a theorem, a subsequence, still denoted by $(\gamma_{\varepsilon_n})$, converges $h$-uniformly to a curve $\gamma=(x,t)$ joining $p$ to $q$.

Choose a finite subdivision of $I$ so fine that, for each of its closed subintervals $[a,b]$, the set $x([a,b])$ lies in the domain $U$ of a coordinate chart of $S$, and choose an open set $V$ with compact closure $\overline V\subset U$ such that $x([a,b])\subset V$. For $n$ large enough, $x_{\varepsilon_n}([a,b])\subset V$ as well, and, in the coordinates of $U\times\R$,
\[
\dot\gamma_{\varepsilon_n}^k(s)
=
\dot\gamma_{\varepsilon_n}^k(a)
-
\int_a^s
\Gamma^k_{\varepsilon_n,ij}(\gamma_{\varepsilon_n}(r))
\dot\gamma_{\varepsilon_n}^i(r)\dot\gamma_{\varepsilon_n}^j(r)\,dr ,\qquad s\in[a,b],
\]
where $\Gamma^k_{\varepsilon_n,ij}$ and $\Gamma^k_{ij}$ denote the Christoffel symbols of $g_{\varepsilon_n}$ and $g$ in these coordinates. Since $g_{\varepsilon_n}=g-\varepsilon_n\,dt^2\to g$ in $C^\infty$ on compact subsets, $\Gamma^k_{\varepsilon_n,ij}\to\Gamma^k_{ij}$ uniformly on $\overline V\times[t_0-C,t_0+C]$; in particular, they are uniformly bounded there. Together with the velocity bounds, this shows that the coordinate velocities $\dot\gamma_{\varepsilon_n}^k$ are uniformly bounded and uniformly Lipschitz on $[a,b]$. By the Ascoli--Arzel\`a theorem, after passing to a further subsequence, common to the finitely many subintervals, we obtain $\gamma_{\varepsilon_n}\to\gamma$ in $C^1(I,S\times\R)$. Passing to the limit in the integral equation above gives
\[
\dot\gamma^k(s)
=
\dot\gamma^k(a)
-
\int_a^s
\Gamma^k_{ij}(\gamma(r))\dot\gamma^i(r)\dot\gamma^j(r)\,dr ,\qquad s\in[a,b].
\]
The right-hand side is continuously differentiable in $s$, so $\gamma$ is $C^2$ on $[a,b]$ and satisfies the geodesic equation of $g$ there. Hence $\gamma$ is a geodesic of $(M,g)$ joining $p$ to $q$.\qedhere
\end{proof}

\section{Consequences and examples}\label{consequences}
We first make explicit a consequence of the proof of Theorem~\ref{main}, which will be used repeatedly in this section.

\begin{corollary}\label{admissibility-equivalence}
Under the hypotheses of Theorem~\ref{main}, fix a complete Killing splitting $M\simeq S\times\R$ and let $p=(x_0,t_0)$, $q=(x_1,t_1)$. Then $p$ and $q$ are joined by a geodesic if and only if there exists an admissible pair $(\sigma,c)$ for $x_0,x_1$.
\end{corollary}
\begin{proof}
If $p$ and $q$ are joined by a geodesic, then $\mathcal C_K(p,q)\ne\emptyset$, and Remark~\ref{admissible_from_Ck} provides an admissible pair. Conversely, the proof of Theorem~\ref{main} uses the hypothesis $\mathcal C_K(p,q)\ne\emptyset$ only through Lemma~\ref{uniform2} based on the existence of an admissible pair.
\end{proof}
Thus, under the hypotheses of Theorem~\ref{main}, the splitting-dependent admissibility condition is equivalent to the intrinsic condition $\mathcal C_K(p,q)\ne\emptyset$.

\begin{remark}\label{variational-criterion}
The admissibility condition has a natural variational meaning. Extend the
definition of $U_{\sigma,c}$ in \eqref{Usigma} to every
$\sigma\in\Omega_{x_0,x_1}(S)$ and $c\in\R$, so that $U_{\sigma,c}=+\infty$ when
$Q_\beta(\sigma,c)=+\infty$. Then
\[
\lim_{\varepsilon\to0^+}J_\varepsilon(\sigma)=\min_{c\in\R}U_{\sigma,c},
\]
both sides being possibly $+\infty$. In particular, $\sigma$ is admissible for
some $c$ if and only if $J_\varepsilon(\sigma)$ stays bounded from above as
$\varepsilon\to0^+$. Let $L\in(-\infty,+\infty]$ be the limit, which exists by Lemma~\ref{monotone}. Since $J_\varepsilon(\sigma)$ increases as $\varepsilon$ decreases, $J_\varepsilon(\sigma)\le L$ for every $\varepsilon\in(0,\varepsilon_*]$. By Lemma~\ref{uniform}, $L\le\inf_cU_{\sigma,c}$. If $L=+\infty$, then $U_{\sigma,c}=+\infty$ for every $c$, and the identity holds trivially. If $L<+\infty$, it suffices to find $\bar c\in\R$ with $U_{\sigma,\bar c}\le L$. For such a $\bar c$, we have $U_{\sigma,\bar c}=L=\inf_cU_{\sigma,c}$, so the infimum is attained. To this end, let $c_\varepsilon$ be the minimizer of
$\mathcal H_\varepsilon(\sigma,\cdot)$ in \eqref{mathcalH}, and choose $b>0$ such that $\beta+\varepsilon_*\le b$
along $\sigma$, which is possible since $\sigma(I)$ is compact. Then
\[
\mathcal H_\varepsilon(\sigma,c)\ge c\Delta t+\frac1{2b}\int_0^1\bigl(\omega(\dot\sigma)-c\bigr)^2\ds,
\]
and the right-hand side is coercive in $c$ and independent of $\varepsilon$. Since
\[
\mathcal H_\varepsilon(\sigma,c_\varepsilon)=J_\varepsilon(\sigma)-\frac12\int_0^1|\dot\sigma|_0^2\ds
\le L-\frac12\int_0^1|\dot\sigma|_0^2\ds,
\]
the family $(c_\varepsilon)$ is bounded,
and $c_{\varepsilon_n}\to\bar c$ along some sequence $\varepsilon_n\to0$. For fixed
$\eta\in(0,\varepsilon_*]$ and $\varepsilon_n\le\eta$ we have $\beta_\eta\ge\beta_{\varepsilon_n}$, hence
$\mathcal H_\eta(\sigma,c_{\varepsilon_n})\le\mathcal H_{\varepsilon_n}(\sigma,c_{\varepsilon_n})
\le L-\frac12\int_0^1|\dot\sigma|_0^2\ds$. Since $\beta_\eta\ge\eta>0$, the map $c\mapsto\mathcal H_\eta(\sigma,c)$ is continuous. Letting $n\to\infty$, we obtain $\frac12\int_0^1|\dot\sigma|_0^2\ds+\mathcal H_{\eta}(\sigma,\bar c)\le L$.  Thus,    passing to the limit as $\eta\to 0^+$, by monotone convergence we get $U_{\sigma,\bar c}\le L$, as required.

The same limit is the supremum of
\[
E_g(\sigma,t):=\frac12\int_0^1\bigl(|\dot\sigma|_0^2+2\omega(\dot\sigma)\dot t-\beta(\sigma)\dot t^2\bigr)\ds
\]
over all $t\in H^1(I,\R)$ with $t(0)=t_0$ and $t(1)=t_1$. Indeed, for fixed $\varepsilon>0$ the map
$t\mapsto E_{g_\varepsilon}(\sigma,t)=E_g(\sigma,t)-\frac\varepsilon2\int_0^1\dot t^2\ds$
is strictly concave on this affine space. Indeed, it is the sum of an affine function of $t$ and of $-\frac12\int_0^1\beta_\varepsilon(\sigma)\dot t^{\,2}\ds$, and the quadratic form $\varphi\mapsto\int_0^1\beta_\varepsilon(\sigma)\dot\varphi^2\ds$ is positive definite on $H^1_0(I,\R)$, since $\beta_\varepsilon\ge\varepsilon$. Its critical points are the functions $t$ such that $\omega(\dot\sigma)-\beta_\varepsilon(\sigma)\dot t$ is constant. Hence its maximum is attained exactly at  $t_\varepsilon$ given by \eqref{Ct}, and \eqref{reazione} gives $\max_tE_{g_\varepsilon}(\sigma,t)=E_{g_\varepsilon}(\sigma,t_\varepsilon)=J_\varepsilon(\sigma)$. Since, for every $t$, $E_{g_\varepsilon}(\sigma,t)$ increases to $E_g(\sigma,t)$ as $\varepsilon\to0^+$, we get
\[
\lim_{\varepsilon\to0^+}J_\varepsilon(\sigma)=\sup_{\varepsilon>0}\sup_tE_{g_\varepsilon}(\sigma,t)=\sup_t\sup_{\varepsilon>0}E_{g_\varepsilon}(\sigma,t)=\sup_tE_g(\sigma,t).
\]
For $\varepsilon=0$ the supremum need not be attained. Indeed, a maximizer $t$ of $E_g(\sigma,\cdot)$ would satisfy $\omega(\dot\sigma)-\beta(\sigma)\dot t\equiv C$ for some $C\in\R$, i.e.\ $(\sigma,t)$ would be an $H^1$ curve with constant Killing charge. For the admissible pair $(\sigma,1)$ of Remark~\ref{no-lift} the supremum is finite, but this equation gives $\dot t=(e^s-C)/(e^s-1)^2$, which is not in $L^2(I)$ for any $C\in\R$.

The proof of Lemma~\ref{uniform2} and the estimates of Section~\ref{foce} use the admissible pair, as well as the minimality of $x_\varepsilon$, only through the bound $J_\varepsilon(x_\varepsilon)\le U_{\sigma,c}$. Two consequences follow. First, in analogy with the characterization of admissible curves above, $p$ and $q$ are joined by a geodesic if and only if $\min J_\varepsilon$ stays bounded from above as $\varepsilon\to0^+$; like admissibility, this criterion is formulated in the chosen splitting, but, being equivalent to~\ref{main-i}, its validity does not depend on that choice. Second, the same estimates apply to any family of critical points of $J_\varepsilon$ whose critical values are bounded from above uniformly in $\varepsilon$; multiplicity results in the spirit of \cite[Theorem~8.4]{BaCaFl17} will be considered elsewhere.
\end{remark}

By reparametrization, the constant Killing charge condition can be weakened to a sign condition.
\begin{corollary}\label{sign-criterion}
Under the hypotheses of Theorem~\ref{main}, two points $p,q\in M$ are joined
by a geodesic if and only if there is a $C^1$ curve $\gamma:I\to M$ joining
them for which $g(\dot\gamma,K)$ is nowhere zero or identically zero.
\end{corollary}
\begin{proof}
A geodesic has constant Killing charge, so it satisfies the stated condition.
Conversely, write $f(s)=g(\dot\gamma(s),K_{\gamma(s)})$.
If $f\equiv0$, then $\gamma\in\mathcal C_K(p,q)$ and Theorem~\ref{main}
applies. If $f$ is nowhere zero, it has constant sign. Set
\[
\Lambda=\int_0^1f(r)\,dr,
\qquad
\psi(s)=\frac1\Lambda\int_0^s f(r)\,dr.
\]
Then $\psi'(s)=f(s)/\Lambda>0$, so $\psi$ is a $C^1$ diffeomorphism of
$I$ onto itself. The curve $\widetilde\gamma=\gamma\circ\psi^{-1}$ has
constant Killing charge $\Lambda$ and joins the same endpoints. Hence
$\widetilde\gamma\in\mathcal C_K(p,q)$, and Theorem~\ref{main} applies again.
\end{proof}

\begin{corollary}\label{orbit-invariance}
Under the hypotheses of Theorem~\ref{main}, let $\Phi$ be the flow of $K$.
For any $p,q\in M$ and any $a,b\in\R$, the points $p,q$ are joined by a
geodesic if and only if $\Phi_a(p),\Phi_b(q)$ are joined by a geodesic.
\end{corollary}
\begin{proof}
In a complete Killing splitting the flow changes only the temporal coordinate, so the pairs $p,q$ and $\Phi_a(p),\Phi_b(q)$ have the same spatial components $x_0,x_1$. By Corollary~\ref{admissibility-equivalence}, in both cases geodesic connectability is equivalent to the existence of an admissible pair for $x_0,x_1$.
\end{proof}
Note that the shifts $a$ and $b$ in Corollary~\ref{orbit-invariance} are independent.
\begin{corollary}\label{sufficient-criteria}
Under the hypotheses of Theorem~\ref{main}, fix a complete Killing splitting
$M\simeq S\times\R$. Let $x_0,x_1\in S$ satisfy either of the following:
\begin{enumerate}[label=\textup{(\alph*)},leftmargin=1.2cm]
\item $x_0,x_1$ belong to the same connected component of $\{\beta>0\}$;
\item there is a curve $\sigma\in\Omega_{x_0,x_1}(S)$ such that
$\omega(\dot\sigma)\equiv c$ almost everywhere for some $c\in\R$.
\end{enumerate}
Then $(x_0,t_0)$ and $(x_1,t_1)$ are joined by a geodesic for every
$t_0,t_1\in\R$.
\end{corollary}
\begin{proof}
In case (a), choose a smooth joining curve $\sigma$ contained in the
component. The compactness of its image gives
$\min_I\beta(\sigma)>0$, hence $Q_\beta(\sigma,0)<+\infty$.
In case~(b), $Q_\beta(\sigma,c)=0$ by definition.
In both cases an admissible pair for $x_0,x_1$ exists, and Corollary~\ref{admissibility-equivalence} gives the conclusion.
\end{proof}
We notice that no completeness assumption on the component of $\{\beta>0\}$ is needed, since
Theorem~\ref{main} is applied to the original spacetime.
In particular, if $\beta>0$ on $S$, condition~(a) holds for every
pair of spatial endpoints, and we recover the result of \cite[Theorem~1.1]{CaFlSa2008} for stationary spacetimes.
Condition~(a) is intrinsic, since $\{g(K,K)<0\}=\{\beta>0\}\times\R$; it means that $(x_0,t_0)$ and $(x_1,t_1)$ lie in the same connected component of the open set where $K$ is timelike. Condition~(b) depends instead on the splitting, since it involves curves in the slice $S$; as $\omega(v)=g(v,K)$ for $v\in TS$, it means that $S$ contains a curve with constant Killing charge joining the $K$-orbits of the two points.

\begin{corollary}\label{nonintegrable}
Under the hypotheses of Theorem~\ref{main}, assume that, in some complete Killing splitting, the distribution $\mathcal D:=\ker\omega$ is nowhere integrable, i.e.\ $(\omega\wedge d\omega)_x\neq0$ for all $x\in S$. Then $(M,g)$ is geodesically connected.
\end{corollary}
\begin{proof}
For local sections $X,Y$ of $\mathcal D$ one has $\omega([X,Y])=-d\omega(X,Y)$, and $(\omega\wedge d\omega)_x\neq0$ implies that $d\omega_x$ does not vanish on $\mathcal D_x\times\mathcal D_x$. Hence $\mathcal D$ is bracket generating and, since $S$ is connected, the Chow--Rashevskii theorem gives, for any $x_0,x_1\in S$, a piecewise smooth curve $\sigma\in\Omega_{x_0,x_1}(S)$ tangent to $\mathcal D$, i.e.\ with $\omega(\dot\sigma)\equiv0$. Corollary~\ref{sufficient-criteria}-(b) gives the conclusion.
\end{proof}
Thus the conclusion of \cite[Theorem~1]{BC} holds without the upper bound on $\beta$, the orthogonality condition $g_0(\nabla\beta,\omega^\sharp)=0$ and the two-sided bound on $|\omega|_0$ assumed there.

In the next two examples $S=\R^3$, $g_0$ is the Euclidean metric, $K=\partial_t$, and $\omega$ and $\beta$ are bounded. Under these conditions, $(M,g)$ is globally hyperbolic and the slices $t=\mathrm{constant}$ are complete spacelike Cauchy hypersurfaces. Indeed, as stated in the Introduction, $t$ is strictly increasing along future-directed causal curves. If such a curve is parametrized by $t$, its spatial velocity $v$ satisfies $|v|_0^2+2\omega(v)-\beta\le0$, hence $|v|_0\le|\omega|_0+\sqrt{|\omega|_0^2+\beta}$, which is bounded. If an inextendible causal curve had a finite upper or lower endpoint for its $t$-range, its spatial component would have a limit there, by this uniform Lipschitz bound, and the curve would have an endpoint in $M$, a contradiction. Hence every inextendible causal curve meets each slice exactly once. Each slice is spacelike and complete, since its induced metric is $g_0$.

\begin{example}\label{positive-example}
Consider $M=\R^3\times\R$, with coordinates $(u,y,z,t)$, and the metric
\begin{equation*}
g=du^2+dy^2+dz^2+2\,du\,dt-\bigl(1-e^{-u^2}\bigr)dt^2.
\end{equation*}
Here $g_0$ is Euclidean, $\omega=du$, and 	$\beta(u)=1-e^{-u^2}$. In particular, 	$\beta+|\omega|_0^2=\beta+1>0$, so $g$ is Lorentzian. 	The vector field $K=\partial_t$ is complete and Killing, 	timelike where $u\ne0$, and lightlike on $\{u=0\}$.
	
Since $\omega$ and $\beta$ are bounded, $(M,g)$ is globally hyperbolic, with complete spacelike Cauchy slices $t=\mathrm{constant}$, as observed above. We show that it is geodesically connected.
	
Fix arbitrary points $p=(u_0,y_0,z_0,t_0)$ and $q=(u_1,y_1,z_1,t_1)$, and let $\Delta t=t_1-t_0$.
Geodesic connectedness follows from Corollary~\ref{sufficient-criteria}-(b). 	Indeed, the affine spatial segment $\sigma$ joining $(u_0,y_0,z_0)$ to $(u_1,y_1,z_1)$ satisfies 	$\omega(\dot\sigma)\equiv u_1-u_0$, hence 	$Q_\beta(\sigma,u_1-u_0)=0$.
	
The admissible pair $(\sigma, u_1-u_0)$ also illustrates Remark~\ref{no-lift}. If $\Delta t\ne0$ and $(u_0,u_1)\ne(0,0)$, then $\sigma$ 	admits no $H^1$ lift $\gamma$ joining $p$ to $q$ with 
$g(\dot\gamma, K)\equiv u_1-u_0$. Indeed, $\beta(u)>0$ along $\sigma$ except possibly 	at one parameter value, so \eqref{splitting-constant} forces  $\dot t=0$ almost everywhere, contradicting $\Delta t\ne0$. Both restrictions are necessary: if $\Delta t=0$, a constant 	time component gives such a lift; if $u_0=u_1=0$, then 	$\beta\equiv0$ along $\sigma$ and any $H^1$ time component with the prescribed endpoints does so. Nevertheless, $\mathcal C_K(p,q)\ne\emptyset$, since it	contains a geodesic.
	
Note that, in this splitting, $\ker du$ is integrable and $g_0(\nabla\beta,\omega^\sharp)=2ue^{-u^2}$. Thus, although the coefficient bounds required in \cite[Theorem~1]{BC} are satisfied, its non-integrability and orthogonality hypotheses are not. Moreover, Corollary~\ref{nonintegrable} does not apply in any complete Killing splitting. Indeed, the one-form $\alpha=g(K,\cdot)=du-\beta(u)\,dt$ satisfies $d\alpha=-\beta'(u)\,du\wedge dt$, hence $\alpha\wedge d\alpha=0$ on $M$. In any complete Killing splitting, with slice $\iota\colon S'\hookrightarrow M$, the corresponding one-form is $\omega'=\iota^*\alpha$, since $\omega'(v)=g(v,K)$ for $v\in TS'$; therefore $\omega'\wedge d\omega'=\iota^*(\alpha\wedge d\alpha)=0$.
\end{example}

\begin{example}\label{negative-example}
Condition~\ref{main-ii} in Theorem~\ref{main} can fail for some pairs 	even when $K$ is timelike on a non-empty open set. For lightlike $K$, an example of this phenomenon is given in \cite[Section~7(c)]{BaCaFl17}.	Choose a smooth non-decreasing function $\chi:\R\to[0,1]$ such that
\[
\chi=0\ \text{on }(-\infty,-1],\qquad	\chi=1\ \text{on }[1,+\infty),
\qquad	0<\chi<1\ \text{on }(-1,1),
\]
and put
\[
\lambda=1-2\chi,\qquad	\beta=4\chi(1-\chi)=1-\lambda^2.\]
On $M=\R^3\times\R$, with coordinates $(x,y,z,t)$, consider
\begin{align*}
g&=dx^2+dy^2+dz^2+2\lambda(x)\,dx\,dt-\beta(x)dt^2\\
&=\bigl(dx+\lambda(x)dt\bigr)^2+dy^2+dz^2-dt^2.
\end{align*}
The $(x,t)$ block has determinant $-1$, so $g$ is Lorentzian. 	The field $K=\partial_t$ is complete and Killing, timelike for 	$-1<x<1$ and lightlike for $x\le-1$ or $x\ge1$. Since $\omega=\lambda\,dx$ and $\beta$ are bounded, $(M,g)$ is globally hyperbolic, with complete spacelike Cauchy slices $t=\mathrm{constant}$, as observed above Example~\ref{positive-example}.
	
Nevertheless, for arbitrary $y_0,z_0,y_1,z_1,t_0,t_1\in\R$, the points
\[
p=(-2,y_0,z_0,t_0),\qquad q=(2,y_1,z_1,t_1)
\]
are not joined by a geodesic. Indeed, suppose that an $H^1$ curve 	$\gamma=(x,y,z,t)$ joining them satisfies \eqref{splitting-constant} for some  $c\in \R$.	Since $x(0)=-2$ and $x(1)=2$, continuity ensures that $a:=\min\{s\in[0,1]:x(s)=-1\}$
is well defined and belongs to $(0,1)$. Moreover, $x(s)<-1$ for every $s<a$. 
Thus $\lambda(x)=1$ and $\beta(x)=0$ on $[0,a]$, so \eqref{splitting-constant} yields $\dot x=c$ almost everywhere there. Since $x$ is absolutely continuous,
$1=x(a)-x(0)=ca$, which implies $c>0$. 	Similarly, let $b\in(0,1)$ be the last hitting time of $1$. On $[b,1]$ we have $x\ge1$, $\lambda=-1$ and $\beta=0$; hence	$-\dot x=c$ almost everywhere and $1=x(1)-x(b)=-c(1-b)$, which implies $c<0$, a contradiction.
Thus $\mathcal C_K(p,q)=\emptyset$, ruling out a connecting geodesic. On the other hand, any two spatial points in the slab 	$(-1,1)\times\R^2=\{\beta>0\}$ give geodesically connectable spacetime endpoints for every choice of their temporal coordinates, 	by Corollary~\ref{sufficient-criteria}-(a). 

We next exhibit a pair of points not covered by either condition of that corollary in the present splitting. Set
\[
x_*=\sup\{x\in\R:\chi(x)\le\tfrac12\}.
\]
Since $\chi$ is continuous and non-decreasing, $x_*\in(-1,1)$, 	$\lambda\ge0$ on $(-\infty,x_*]$, and $\lambda<0$ on 	$(x_*,+\infty)$. 	Fix $\bar x\in(x_*,1)$ and let
$r:=(\bar x,\bar y,\bar z,\bar t)$, with $\bar y,\bar z,\bar t\in\R$ arbitrary.
Let $\sigma=(x,y,z)$ be the affine spatial segment from 	$(-2,y_0,z_0)$ to $(\bar x,\bar y,\bar z)$, so that
\[
x(s)=-2+cs,\qquad c=\bar x+2>0.
\]
On $\{s:x(s)>-1\}$ set
\[
\dot t=-\frac{c}{2(1-\chi(x))},
\]
which is bounded because $\chi(x)\le\chi(\bar x)<1$ along the segment; 	a direct computation gives	$\lambda(x)\dot x-\beta(x)\dot t=c$ there. 	On the initial interval $\{s:x(s)\le-1\}=[0,1/c]$ we have 	$\lambda=1$ and $\beta=0$, so the constraint reduces to 	$\dot x=c$, which holds, while $\dot t$ is unrestricted. 	Choosing $\dot t$ equal to a suitable constant there gives a 	Lipschitz curve $\gamma=(\sigma,t)$ from $p$ to $r$ with 	Killing charge $c$. 	Hence $\gamma\in\mathcal C_K(p,r)$, and Theorem~\ref{main} 	shows that $p$ and $r$ are joined by a geodesic, for every 	choice of $t_0$ and $\bar t$.

The pair $(\sigma,c)$ is admissible for the spatial components of $p$ and $r$, and
\[
0<Q_\beta(\sigma,c)	=\int_0^1\beta(x)\dot t^{\,2}\ds	=c\int_{-1}^{\bar x}\frac{\chi(x)}{1-\chi(x)}\,dx<+\infty.
\]
Neither condition of Corollary~\ref{sufficient-criteria} applies to these endpoints. Condition~(a) fails because
$\beta(-2)=0$. To exclude condition~(b), let $\Lambda$ 	be a primitive of $\lambda$. 	If $\widetilde\sigma=(\widetilde x,\widetilde y,\widetilde z) 	\in H^1(I,S)$ joins the two spatial endpoints, then
$\frac{d}{ds}\Lambda(\widetilde x)=\omega(\dot{\widetilde\sigma})$ almost everywhere. Moreover $\widetilde x$ takes the value $x_*$. If 	$\omega(\dot{\widetilde\sigma})\equiv\widetilde c$, 	then  $\Lambda\circ\widetilde x$ is  affine. On the other hand, the sign of $\lambda$, together with $\lambda=1$ on $[-2,-1]$, gives
\[
\Lambda(x_*)>\max\{\Lambda(-2),\Lambda(\bar x)\},
\] 	
in contradiction with  $\Lambda\circ\widetilde x$ being affine.

Since $Q_\beta(\widetilde\sigma,\widetilde c)=0$ exactly when $\omega(\dot{\widetilde\sigma})=\widetilde c$ almost everywhere, these spatial endpoints admit admissible pairs, but none with $Q_\beta=0$. Thus, in the chosen splitting, the integral admissibility condition is strictly more general than the sufficient conditions in Corollary~\ref{sufficient-criteria}. Condition~(b), however, may hold for the same endpoints in another complete Killing splitting.

Similarly, $r$ and $q$ are joined by a curve with Killing charge 	$-\delta$, where $\delta=2-\bar x$. 	Take the affine spatial segment from $(\bar x,\bar y,\bar z)$ to 	$(2,y_1,z_1)$, so that $x(s)=\bar x+\delta s$, and set
\[
\dot t=\frac{\delta}{2\chi(x)}\qquad\text{where }x<1.
\]
This derivative is bounded because $\chi(x)\ge\chi(\bar x)>\tfrac12$. 	Choose $\dot t$ constant on the final interval where $x\ge1$ 	to match the temporal endpoints. 	Unlike the first one, this second pair is covered by Corollary~\ref{sufficient-criteria}-(b), since 	$\lambda<0$ on $[\bar x,2]$. 	
Consequently, $p$ and $q$ are both joined to $r$ by geodesics, 	although they are not joined to each other. The obstruction comes from the opposite signs of the Killing charges, as explained below. 
\end{example}
\begin{remark}
Admissible pairs whose constants are non-zero and have the same sign can 	be concatenated. Let $(\sigma_1,c_1)$ and $(\sigma_2,c_2)$ be 	admissible for their respective endpoints, with 	$\sigma_1(1)=\sigma_2(0)$ and $c_1c_2>0$. Set 
\[
\theta=\frac{c_1}{c_1+c_2}\in(0,1),\qquad C=c_1+c_2,
\]
and define
\[
\sigma(s)=	\begin{cases}
		\sigma_1(s/\theta),&0\le s\le\theta,\\
		\sigma_2((s-\theta)/(1-\theta)),&\theta\le s\le1.
	\end{cases}
\]
Since $\mathfrak q(x,kv,kc)=k^2\mathfrak q(x,v,c)$ for $k>0$,
\[
Q_\beta(\sigma,C)=\frac{Q_\beta(\sigma_1,c_1)}{\theta}	+\frac{Q_\beta(\sigma_2,c_2)}{1-\theta}<+\infty,
\]
so $(\sigma,C)$ is admissible for the endpoints 	$\sigma_1(0)$ and $\sigma_2(1)$.
If $c_1=c_2=0$, the same construction works with any 	$\theta\in(0,1)$ and $C=0$.
Consequently, under the hypotheses of Theorem~\ref{main}, Remark~\ref{admissible_from_Ck} and Corollary~\ref{admissibility-equivalence} 	show that if $p$ is joined to $r$ and $r$ to $q$ by curves with
non-zero constant Killing charges of the same sign, or both zero, 	then $p$ and $q$ are joined by a geodesic.
	
For the points $p$, $r$, $q$ of Example~\ref{negative-example}, the arguments based on the first and last hitting times show that every $H^1$ curve with constant Killing charge
from $p$ to $r$ has positive charge, whereas every such curve from $r$ to $q$ has negative charge. Any curve from  $p$ to $q$ with constant Killing charge would have to satisfy both sign requirements, regardless of whether it passes through $r$. This incompatibility rules out a geodesic joining $p$ to $q$.
\end{remark}

\section*{Acknowledgments}
\noindent We thank Miguel S\'anchez for his comments on a preliminary version of this paper and, in particular, for suggesting the variational interpretation of the admissibility condition in Remark~\ref{variational-criterion}. Both authors  are  partially supported   by  MUR under the Programme ``Department of Excellence'' Legge 232/2016  (Grant No. CUP - D93C23000100001) and by
``INdAM - GNAMPA Project''  CUP E53C25002010001.



\begin{thebibliography}{99}



\bibitem{AtGuHo20}
C.~Atindogbe, M.~Guti\'errez and R.~Hounnonkpe,
\emph{Lorentzian manifolds with causal Killing vector field: causality and geodesic connectedness},
Ann. Mat. Pura Appl. \textbf{199} (2020), 1895--1908.

\bibitem{BaCaFl17}
R.~Bartolo, A.~M.~Candela and J.~L.~Flores,
\emph{Connectivity by geodesics on globally hyperbolic spacetimes with a lightlike Killing vector field},
Rev. Mat. Iberoam. \textbf{33} (2017),  1--28.



\bibitem{BC}
R.~Bartolo and E.~Caponio,
\emph{Trajectories of affine control systems and geodesics of a spacetime with a causal Killing vector field},
J. Dyn. Control Syst. \textbf{29} (2023), 1409--1424.

\bibitem{BC24}
R.~Bartolo and E.~Caponio,
\emph{Correction to: Trajectories of affine control systems and geodesics of a spacetime with a causal Killing vector field},
J. Dyn. Control Syst. \textbf{30} (2024), Paper No.~45.


\bibitem{BeEhEa96}
{ J.~K. Beem, P.~E. Ehrlich and K.~L. Easley}, {\em {Global {L}orentzian
		Geometry}},  Marcel Dekker Inc., New York, second~ed., 1996.

\bibitem{BeFoGi91}
V.~Benci, D.~Fortunato and F.~Giannoni,
\emph{On the existence of multiple geodesics in static space-times},
Ann. Inst. H. Poincar\'e Anal. Non Lin\'eaire \textbf{8} (1991), 79--102.

\bibitem{BernalSanchez03}
A.~N.~Bernal and M.~S\'anchez,
\emph{On smooth Cauchy hypersurfaces and Geroch's splitting theorem},
Comm. Math. Phys. \textbf{243} (2003),  461--470.



\bibitem{CaFlSa2008}
A.~M.~Candela, J.~L.~Flores and M.~S\'anchez,
\emph{Global hyperbolicity and Palais--Smale condition for action functionals in stationary spacetimes},
Adv. Math. \textbf{218} (2008),  515--536.

\bibitem{CapCor23}
E.~Caponio and D.~Corona,
\emph{A variational setting for an indefinite Lagrangian with an affine Noether charge},
Calc. Var. Partial Differential Equations \textbf{62} (2023), Paper No.~39.

\bibitem{CapCorGiaPic24}
E.~Caponio, D.~Corona, R.~Giamb\`o and P.~Piccione,
\emph{Fixed energy solutions to the Euler--Lagrange equations of an indefinite Lagrangian with affine Noether charge}, Ann. Mat. Pura Appl. \textbf{203} (2024), 1819--1850.


	\bibitem{CaJaSa14}
{ E.~Caponio, M.~A. Javaloyes, and M.~S\'anchez}, {\em Wind Finslerian
	structures: from Zermelo's navigation to the causality of spacetimes},
Mem. Amer. Math. Soc. \textbf{300} (2024), no.~1501.


\bibitem{Dieudonne70}
J.~Dieudonn\'e,
\emph{Sur un th\'eor\`eme de Glaeser},
J. Analyse Math. \textbf{23} (1970), 85--88.

\bibitem{GiaMas91}
F.~Giannoni and A.~Masiello,
\emph{On the existence of geodesics on stationary Lorentz manifolds with convex boundary},
J. Funct. Anal. \textbf{101} (1991),  340--369.

\bibitem{GiaPic99}
F.~Giannoni and P.~Piccione,
\emph{An intrinsic approach to the geodesical connectedness of stationary Lorentzian manifolds},
Comm. Anal. Geom. \textbf{7} (1999),  157--197.

\bibitem{Glaeser63}
G.~Glaeser,
\emph{Racine carr\'ee d'une fonction diff\'erentiable},
Ann. Inst. Fourier (Grenoble) \textbf{13} (1963),  203--210.


\bibitem{Minguzzi08}
E.~Minguzzi,
\emph{Limit curve theorems in Lorentzian geometry},
J. Math. Phys. \textbf{49} (2008), 092501, 18 pp.

\end{thebibliography}
\end{document}